\documentclass[a4paper,twoside]{IEEEtran}

\usepackage{graphicx}  

\usepackage{url}       

\usepackage{amsmath,amssymb}   
\usepackage[ruled,vlined,lined,linesnumbered,commentsnumbered]{algorithm2e}

\usepackage{color}
\def\Qset{\mathcal{Q}}

\def\PQ{P_{\Qset}}

\def\cn{C_\mathrm{N}}
\def\cgp{C_\mathrm{GP}}

\def\algname{{\bf UPN}}

\def\algnameBB{{\bf GPBB}}
\def\algnameGP{{\bf GP}}

\def\algnameBT{{\bf BT}}

\def\D1{D_1}		

\def\gridsize{64}
\def\xtrue{x_\mathrm{true}}

\newcommand{\xtv}{x^*}
\newcommand{\half}{\frac{1}{2}}
\newcommand{\argmin}{\mathop{\rm argmin}}

\newcommand{\TV}[1]{\| #1 \|_\mathrm{TV}}
\begin{document}
%
\title{Accelerated gradient methods for total-variation-based CT image reconstruction}
%
%
\author{Jakob~H.~J{\o}rgensen, \thanks{Jakob H. J{\o}rgensen and Per Christian Hansen are with the Department of Informa\-tics and Mathematical Modelling, Technical University of Denmark, Richard Petersens Plads, Building 321, 2800 Kgs. Lyngby, Denmark.
Corresponding author: Jakob H. J{\o}rgensen, E-mail:
jakj@imm.dtu.dk.}%
        Tobias~L.~Jensen, \thanks{Tobias L. Jensen and S{\o}ren H. Jensen are with the Department of Electronic Systems, Aalborg University, Niels Jernesvej 12, 9220 Aalborg {\O}, Denmark.}%
        Per~Christian~Hansen,
        S{\o}ren H.~Jensen,
        Emil Y. Sidky,  \thanks{Emil Y. Sidky and Xiaochuan Pan are with the Department of Radiology, The University of Chicago, 5841 S. Maryland Avenue, Chicago, IL 60637, USA.}%
        and Xiaochuan Pan}

%
%
%

\maketitle

\thispagestyle{empty}

\begin{abstract}
Total-variation (TV)-based CT image reconstruction has shown experimentally to be capable of producing accurate reconstructions from sparse-view data.
In particular TV-based reconstruction is very well suited for images with piecewise nearly constant regions.
Computationally, however, TV-based reconstruction is much more demanding, especially for 3D imaging, and the reconstruction from clinical data sets is far from being close to real-time. This is undesirable from a clinical perspective, and thus there is an incentive to accelerate the solution of the underlying optimization problem.

%
The TV reconstruction can in principle be found by any optimization method, but in practice the large scale of the systems arising in CT image reconstruction preclude the use of memory-demanding methods such as Newton's method. The simple gradient method has much lower memory requirements, but exhibits slow convergence.

In the present work we address the question of how to reduce the number of gradient method iterations needed to achieve a high-accuracy TV reconstruction. We consider the use of two accelerated gradient-based methods, \algnameBB{} and \algname{}, to solve the
3D-TV minimization problem in CT image reconstruction. The former incorporates several heuristics from the optimization literature such as Barzilai-Borwein (BB) step size selection and nonmonotone line search. The latter uses a cleverly chosen sequence of auxiliary points to achieve a better convergence rate.

The methods
are memory efficient and equipped with a stopping criterion to ensure that the TV reconstruction has indeed been found.
An implementation of the methods (in C with interface to Matlab) is available for download from  \url{http://www2.imm.dtu.dk/~pch/TVReg/}.
%


We compare the proposed methods with the standard gradient method, applied to a 3D test problem with synthetic few-view data. We find experimentally that for realistic parameters the proposed methods significantly outperform the gradient method. 

\end{abstract}


\section{Introduction}
Algorithm development for image reconstruction from incomplete data has experienced renewed interest in the past years. Incomplete data arises for instance in case of a small number of views, and the development of algorithms for incomplete data thus has the potential for a reduction in imaging time and the delivered dosage.

Total-varation (TV)-based image reconstruction is a promising direction, as experiments have documented the potential for accurate image reconstruction under conditions such as few-view and limited-angle data, see e.g. \cite{Sidky:2006}.

However, it also known that it is difficult to design fast algorithms for obtaining exact TV reconstructions due to nonlinearity and nonsmoothness of the underlying optimization problem.
Many different approaches have been developed, such as
time marching \cite{Rudin92}, fixed-point iteration \cite{Vogel:96},
and various minimization-based methods such as sub-gradient methods
\cite{Combettes:02}, second-order cone programming (SOCP)
\cite{Goldfarb:05}, and duality-based methods
\cite{Chambolle:04,Hintermuller:06} -- but for large-scale applications such as CT image reconstruction the computational burden is still unacceptable. As a consequence heuristic and much faster techniques such as the one in \cite{Sidky:2006} for approximating the TV solution have been developed. 
In such inaccurate, but efficient, TV-minimization solvers the resulting image depends on several algorithm parameters, which introduces an unavoidable variability.
In contrast, for the accurate TV algorithms considered here, the resulting image can be considered dependent only on the parameters of the optimization problem.

In this work we present two accelerated gradient-based optimization methods that are capable of computing the TV reconstruction of 3D volumes to within an accuracy specified by the user.

\section{Theory}

\subsection{Total-variation-based image reconstruction}
In this work we consider total-variation (TV)-based image reconstruction for computed tomography.
The 3D reconstruction is represented by the vector $\xtv$ which is the solution to the minimization problem
  \begin{equation}
  \label{eq:TVobjectivefunction}
    \xtv = \argmin_{x\in \Qset} \ \phi(x), \quad
    \phi(x) = \half \| A\, x-b \|_2^2 + \alpha \TV{x}.
  \end{equation}
Here, $x$ is the unknown image, $\Qset$ is the set of feasible $x$, $A$ is the system matrix, $b$ is the projection data stacked into a column vector, and $\alpha$ is the TV regularization parameter specifying the relative weighting between the fidelity term and the TV term. $\TV{x}$ is the discrete total-variation of $x$,
\begin{equation}
  \label{eq:TV}
    \TV{x} = \sum_{j=1}^N \|D_jx\|_2,
  \end{equation}
where $N$ is the number of voxels and $D_j$ is the forward difference approximation to the gradient at voxel $j$. 

\subsection{Smooth and strongly convex functions}
We recall that a continuously differentiable function $f$ is convex if
  \begin{align}
    f(x) \geq f(y) + \nabla f(y)^T (x-y)
  \end{align}
for all $x$, $y \in \Qset$. A stronger notion of convexity is \emph{strong convexity}: $f$ is said to be
\emph{strongly convex} with \emph{strong convexity parameter} $\mu$ if there exists a  $\mu \geq 0$ such that
  \begin{equation}
  \label{eq:defstronglyconvex}
    f(x) \geq f(y) + \nabla f(y)^T (x-y) + \half\mu \| x-y\|_2^2
  \end{equation}
for all $x$, $y \in \Qset$.
Furthermore, $f$ has \emph{Lipschitz continuous}
gradient with \emph{Lip\-schitz constant} $L$, if
  \begin{equation}
  \label{eq:defsmooth}
    f(x) \leq f(y) + \nabla f(y)^T(x-y) + \half L \|x-y\|_2^2
  \end{equation}
for all $x$, $y \in \Qset$. The ratio $\mu/L$ is important for the convergence rate of gradient methods we will consider.

The problem \eqref{eq:TVobjectivefunction} can be shown \cite{TVReg} to be strongly convex and have Lipschitz continuous gradient in the case where $A$ specifies
a full-rank overdetermined linear system.
In the
rank deficient or underdetermined case, which occurs for instance for few-view data, the strong convexity assumption is violated. However, as we shall see, this turns out not to pose a problem for the gradient methods we consider.

%

\section{Algorithms}

\subsection{Gradient projection methods}
The optimization problem \eqref{eq:TVobjectivefunction} can, in principle, be solved by use of a simple gradient projection (\algnameGP) method
  \begin{align} \label{eq:gradientmethod}
    x^{(k+1)} = \PQ \left(x^{(k)} - \theta_k \nabla f(x^{(k)}) \right ) ,
    \quad k = 0,1,2,\ldots,
  \end{align}
where $\PQ$ denotes projection onto the set $\Qset$ of feasible $x$, and $\theta_k$ is the step size at the
$k$th step. The worst-case convergence rate of \algnameGP{} with $\mu > 0$ and constant step size is
\begin{align}
 f(x^{(k)}) - f^\star \leq \left(1-\frac{\mu}{L}\right)^k \cdot \cgp, \label{eq:GP_convergence_rate}
\end{align}
where $\cgp$ is a constant \cite[\S 7.1.4]{Nemirovsky:1983}.

 For large-scale imaging modalities, such as CT, this slow convergence renders the simple gradient method impractical. On the other hand the simplicity and the low memory requirements of the gradient method remain attractive. Various modifications have been suggested in the optimization literature. For instance, a significant acceleration is often observed empirically if the gradient method is equipped with a Barzilai-Borwein (BB) step size strategy and a nonmonotone line search \cite{BB:88, Birgin:00,Zhu:08,Grippo:86,Raydan:97}, see Algorithm \ref{algo:BB}: \algnameBB{} for a pseudo-code. Empirically we have found $K=2$ and $\sigma=0.1$ to be satisfactory parameter choices. However, it remains unproven that \algnameBB{} achieves a better worst-case convergence rate than \eqref{eq:GP_convergence_rate}.

\begin{algorithm}[tb]\normalsize
\SetKwInOut{Input}{input}\SetKwInOut{Output}{output}\SetKw{Abort}{abort}
\Input{$x^{(0)}$, $K$}
\Output{$x^{(k+1)}$}
$\theta_0 = 1$ \;
\For{$k = 0, 1, 2, \dots$}{
\tcp{BB strategy}
\If{$k>0$}{
$\theta_k \leftarrow \frac{\|x^{(k)} - x^{(k-1)} \|_2^2}{\langle \,
  x^{(k)} - x^{(k-1)}, \nabla f(x^{(k)}) - \nabla f(x^{(k-1)}) \, \rangle}$ \;
  } 
  $\beta\leftarrow 0.95$ \;
  $\bar{x} \leftarrow \PQ( x^{(k)} - \beta \theta_k \nabla f(x^{(k)}) )$ \;
  $\hat{f} \leftarrow \max\{ f(x^{(k)}), f(x^{(k-1)}),\ldots,f(x^{(k-K)}) \}$ \;
  \While{$f(\bar{x}) \geq \hat{f} - \sigma \,
  \nabla f(x^{(k)})^T(x^{(k)}-\bar{x})$}{
  $\beta \leftarrow \beta^2$ \;
  $\bar{x} \leftarrow \PQ(x^{(k)} - \beta\theta_k \nabla f(x^{(k)}))$ \;
  } 
  $x^{(k+1)} \leftarrow \bar{x}$ \; } 
\caption{\algnameBB{}}\label{algo:BB}
\end{algorithm}

\subsection{Nesterov's optimal method} \label{subsec:nesterov}

Nesterov \cite{Ne:04} proposed a gradient-based method that for given $\mu > 0$ achieves the convergence rate
\begin{align}
 f(x^{(k)}) - f^\star \leq \left(1-\sqrt{\frac{\mu}{L}}\right)^k \cdot \cn, \label{eq:Nesterov_convergence_rate}
\end{align}
where $\cn$ is a constant, and he proved the method to be optimal, i.e., that no gradient-based method can achieve better worst-case convergence rate on the class of strongly convex problems. 

Comparing \eqref{eq:GP_convergence_rate} and \eqref{eq:Nesterov_convergence_rate}, we see how the ratio $\mu/L$ affects the predicted worst-case convergence rates: When $\mu/L$ decreases, both convergence rates become slower, but less in \eqref{eq:Nesterov_convergence_rate} due to  the square root. We therefore expect Nesterov's method to show better convergence for smaller $\mu/L$. Small $\mu/L$ arise for instance when the number of views is small, see \cite{TVReg}. 

Nesterov's method requires that both $\mu$ and $L$ are given by the user, and
in order for the method to be convergent $\mu$ must be chosen sufficiently small and $L$ sufficiently large. For real world applications such as CT, $\mu$ and $L$ are seldom known, which makes the method impractical. Taking overly conservative estimates can depreciate the better convergence rate \eqref{eq:Nesterov_convergence_rate}; hence, accurate estimates of $\mu$ and $L$ are important.

\subsection{Estimating $\mu$ and $L$}

A sufficiently large $L$ can be chosen using \emph{back-tracking line search} \cite{beck:2009,Va:09}, see Algorithm \ref{algo:BT}: \algnameBT{} for pseudo-code.
\begin{algorithm}[tb]\normalsize
\SetKwInOut{Input}{input}\SetKwInOut{Output}{output}
\Input{$y,\bar{L}$}
\Output{$x,\tilde L$}
$\tilde L \leftarrow \bar{L}$ \;
$x \leftarrow \PQ \left( y-\tilde L^{-1} \nabla f(y) \right)$ \;
\While{$f(x)> f(y)+\nabla f(y)^T(x-y) + \half \tilde L \|x-y\|_2^2$}{
	$\tilde  L \leftarrow \rho_L \tilde L$ \;
	$x \leftarrow \PQ \left( y-\tilde L^{-1} \nabla f(y) \right)$ \;
} 
\caption{\algnameBT}\label{algo:BT}
\end{algorithm}
Essentially, an estimate $\bar L$ of $L$ is increased by multiplication with a constant $\rho_L>1$ until \eqref{eq:defsmooth} is satisfied.

Accurately estimating $\mu$, such that \eqref{eq:defstronglyconvex} is satisfied globally, is more difficult. Here, we propose a simple and computationally inexpensive heuristic: In the $k$th iteration choose an estimate $\mu_k$ as the largest value of $\mu$ that satisfies \eqref{eq:defstronglyconvex} between $x^{(k)}$ and $y^{(k)}$,
and make the $\mu_k$-sequence non-increasing:
\begin{align}
    \mu_k = \min \left\{\mu_{k-1},\frac{f(x)-f(y)-\nabla f(y)^T(x-y)}{\frac{1}{2}\|x-y
      \|_2^2} \right\} .   \label{eq:mu_k_estimate}
\end{align}

We call the Nesterov method equipped with estimation of $\mu$ and $L$ Unknown Parameter Nesterov (UPN) and pseudo-code is given in Algorithm \ref{algo:UPN}: \algname{}.
%
%

Unfortunately, convergence of \algname{} is not guaranteed, since the estimate \eqref{eq:mu_k_estimate} can be too large. However, we have found empirically that an estimate sufficient for convergence is typically effectively determined by \eqref{eq:mu_k_estimate}.

It is possible to ensure convergence by introducing a restarting procedure  \cite{TVReg} at the price of lowering the convergence rate bound and thereby losing optimality of the method. However, we have found empirically that the restarting procedure is seldom needed, and for realistic parameters the simple heuristic \eqref{eq:mu_k_estimate} is sufficient.

\begin{algorithm}[tb]\normalsize
\SetKwInOut{Input}{input}\SetKwInOut{Output}{output}\SetKw{Abort}{abort}
\Input{$x^{(0)},\bar \mu,\bar L$}
\Output{$x^{(k+1)}$}
$[ x^{(1)},L_{0} ] \leftarrow \algnameBT{}(x^{(0)},\bar{L})$ \label{algoline:upn_bt1}\;
$\mu_0=\bar \mu, \quad y^{(1)} \leftarrow x^{(1)} , \quad \theta_1
\leftarrow \sqrt{\mu_0/L_0}$ \;
\For{$k = 1, 2, \dots$}{
	$[ x^{(k+1)},L_k ] \leftarrow \algnameBT{}(y^{(k)},L_{k-1})$ \;
    $\mu_k \leftarrow \min \bigl\{ \mu_{k-1}, M(x^{(k)},y^{(k)}) \bigr\}$ \;
    $\theta_{k+1} \leftarrow \mathrm{positive}~\mathrm{root}~\mathrm{of}~
    \theta^2 = (1-\theta)\theta_k^2 + (\mu_k/L_k)\, \theta$ \;
    $\beta_k \leftarrow \theta_k(1-\theta_k)/(\theta_k^2 + \theta_{k+1})$ \;
    $y^{(k+1)} \leftarrow x^{(k+1)} + \beta_k (x^{(k+1)} - x^{(k)})$ \;
} 
\caption{\algname{}}\label{algo:UPN}
\end{algorithm}

\subsection{Stopping criterion}
For an unconstrained convex optimization problem such as \eqref{eq:TVobjectivefunction} the norm of the gradient is a measure of closeness to the minimizer through the first-order optimality conditions \cite{BoVa:04}. For a constrained convex optimization problem it is possible to express a similar optimality condition, namely in terms of the \emph{gradient map} defined by
\begin{align}
 G_\nu(x) = \nu \left(x - \PQ\left(x-\nu^{-1}\nabla f(x)\right)\right),
\end{align}
where $\nu>0$ is a scalar. A point $x^\star$ is optimal if and only if $G_\nu(x^\star)=0$ for any $\nu > 0$ \cite{Va:09}. We can use this to design a stopping criterion: Stop the algorithm after iteration $k$ if $\|G_\nu(x^{(k)})\|_2/N \leq \epsilon$, where $\epsilon$ is a user-specified tolerance.

For an under-determined problem, e.g. in the few-view case, the objective function in \eqref{eq:TVobjectivefunction} is nearly flat at the minimizer, which makes it difficult to determine when a sufficiently accurate reconstruction has been found. Here, the gradient map provides a simple, yet sensitive, stopping criterion.

\section{Simulation results and discussion}

\subsection{Simulation setup}
At this point we emphasize that our objective is to obtain an accurate TV reconstruction while reducing the required number of gradient method iterations. We include two reconstructions merely to demonstrate that the methods indeed are successful in solving \eqref{eq:TVobjectivefunction}, thereby reconstructing the desired image. In \cite{TVReg} dependence of the convergence with respect to parameter variation is explored.

To demonstrate and compare the convergence of \algnameGP{}, \algnameBB{} and \algname{} we set up a simple test problem. As test image $\xtrue$ we use 
the threedimensional FORBILD head phantom
discretized into 
$\gridsize^3$ voxels. 
%
%
We simulate a parallel beam geometry with view directions evenly distributed over the unit sphere. Projections are computed as the forward mapping of the discretized image subject to additive Gaussian white \hbox{noise $e$} of relative magnitude
$\|e\|_2/\|A\xtrue\|_2=0.01$, i.e.,
\begin{align}
 b = A\xtrue + e.
\end{align}
We enforce nonnegativity by taking $\Qset = \mathbb{R}_+^{\gridsize^3}$.
We consider two reconstructions: A ``many-view'' using 
55 views and a ``few-view'' using only 
19 views of size 
$91^2$ pixels. In the latter case $A$ has less rows than columns, which can be shown \cite{TVReg} to lead to violation of the assumption on strong convexity. The iterations are continued to the tolerance $\epsilon = 10^{-8}$ is met.

\subsection{Simulation results}
\begin{figure*}
\centerline{
\includegraphics[width=0.24\textwidth]{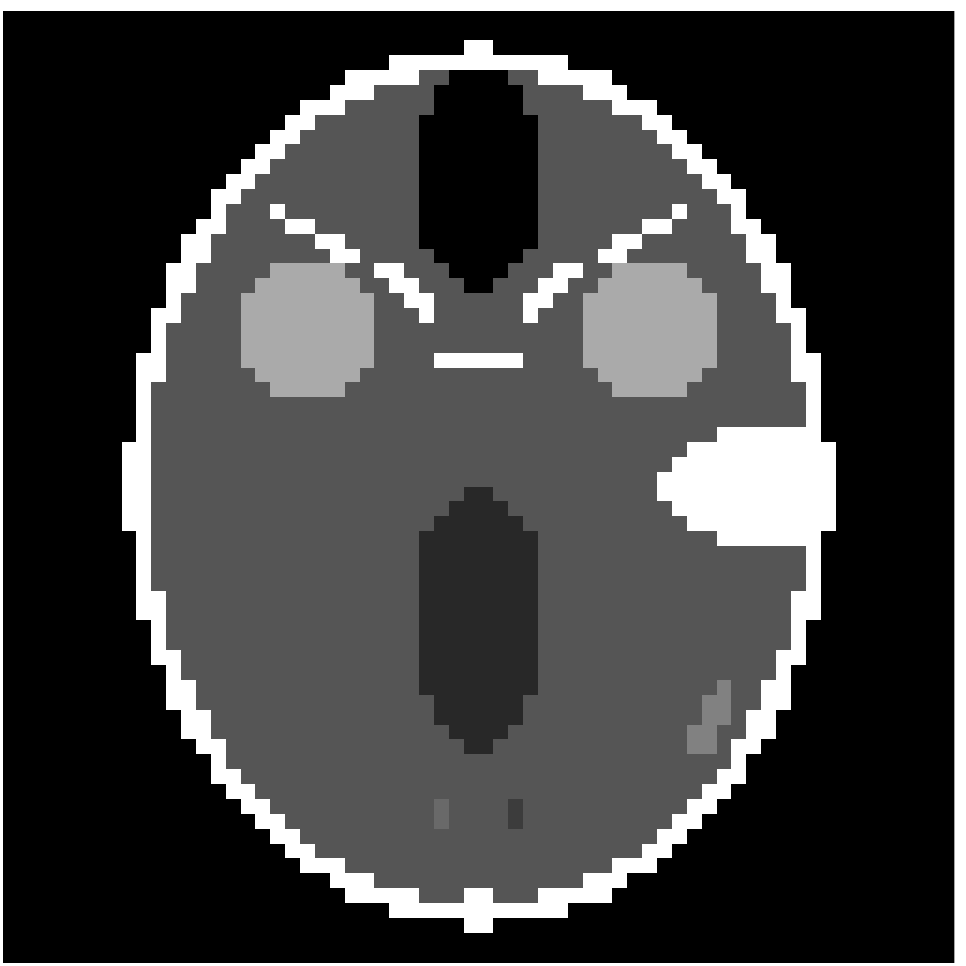}
\includegraphics[width=0.24\textwidth]{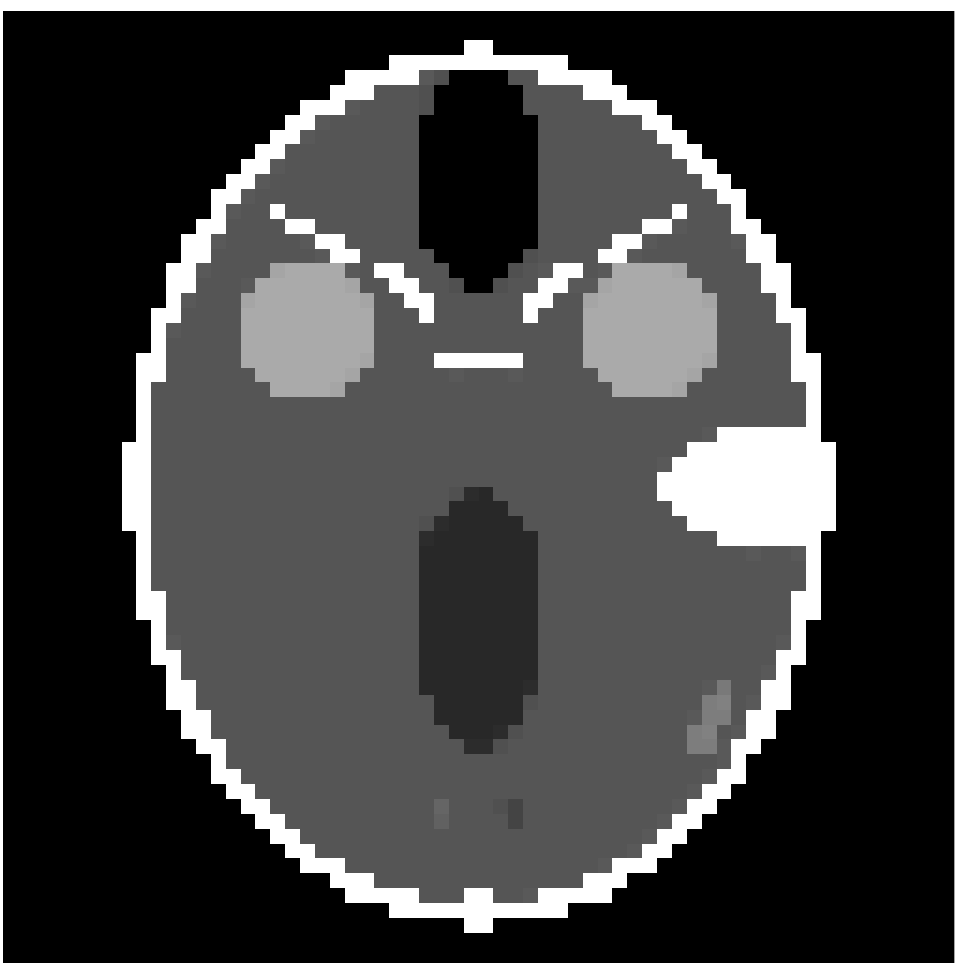}
\includegraphics[width=0.24\textwidth]{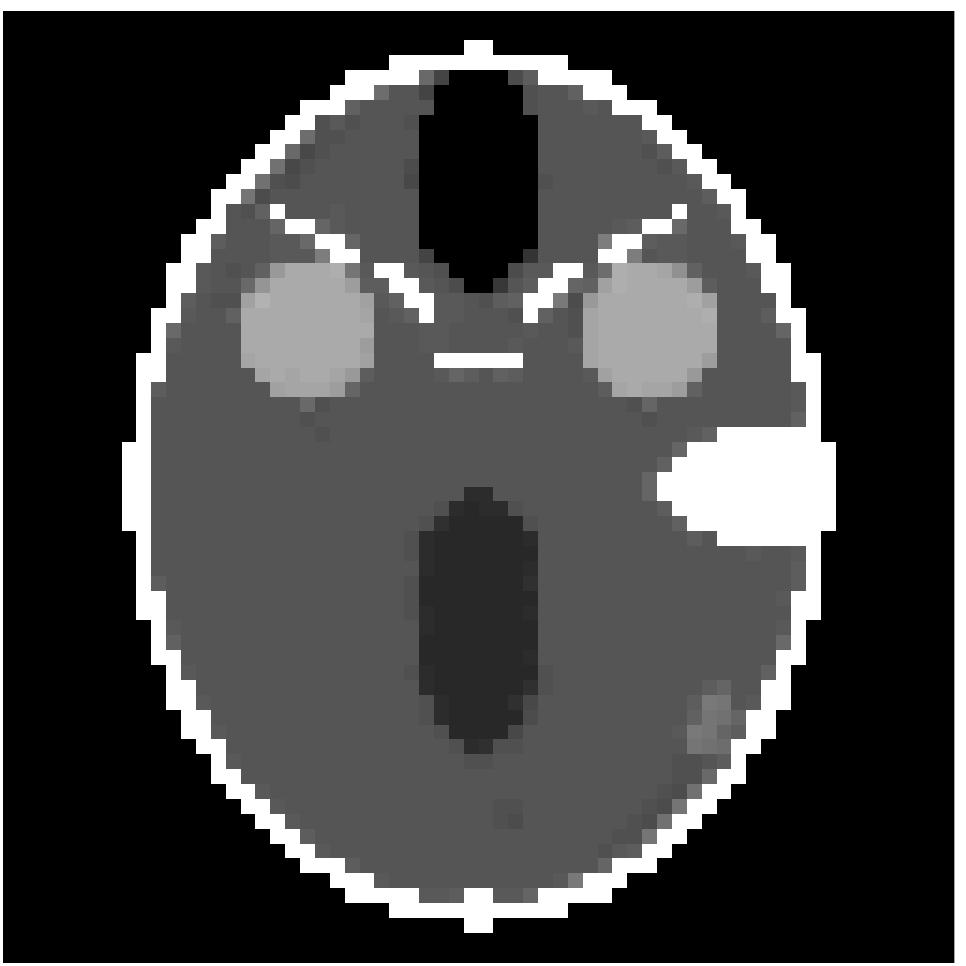}}
\caption{Central axial slices. Left: Original. Center: Many-view \algname{} reconstruction. Right: Few-view \algname{} reconstruction. The display color range is set to $[1.04,1.07]$ for improved viewing contrast.}
\label{fig:reconstructions}
\end{figure*}

\begin{figure*}
\centerline{
\includegraphics[width=0.32\textwidth]{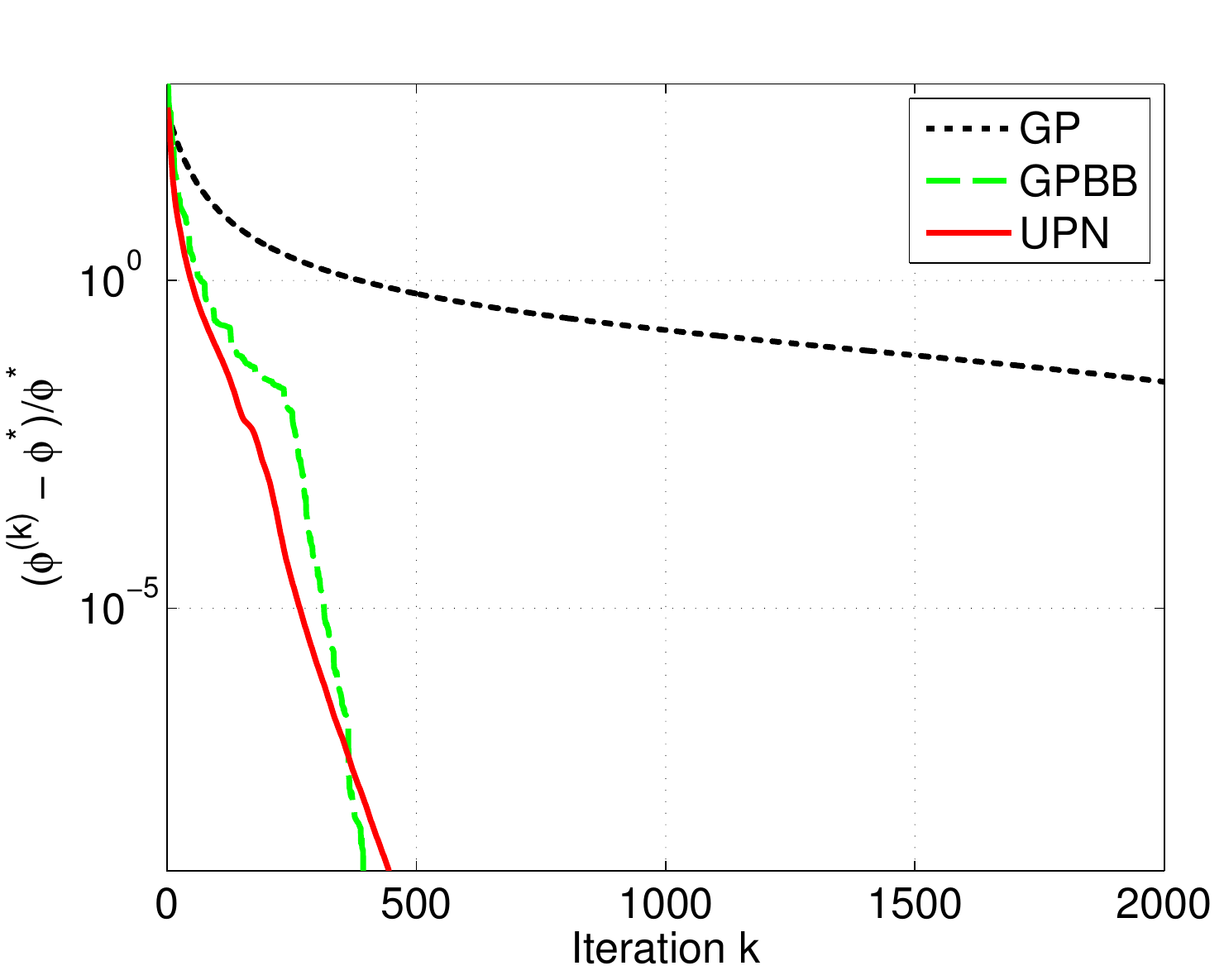}
\includegraphics[width=0.32\textwidth]{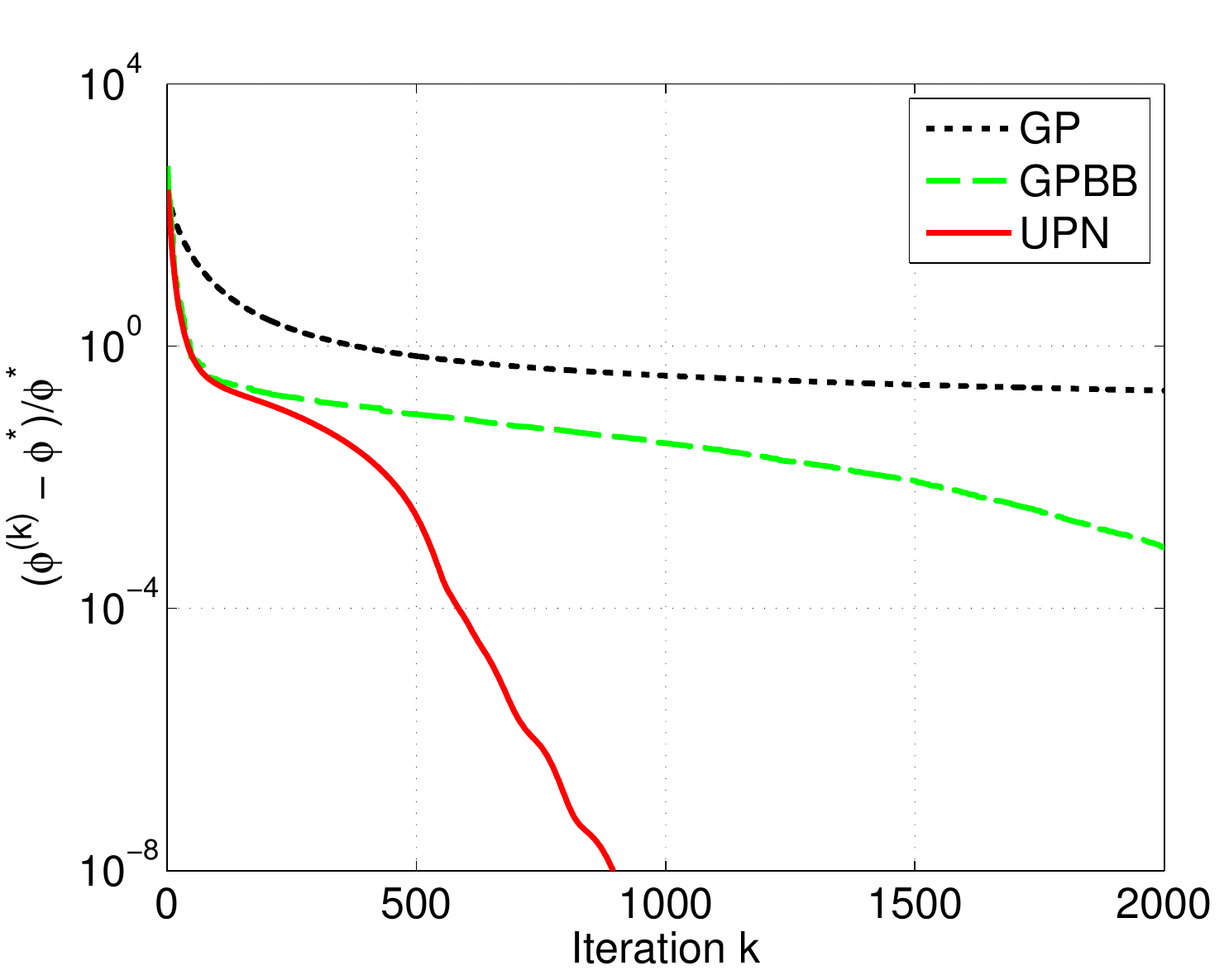}
\includegraphics[width=0.32\textwidth]{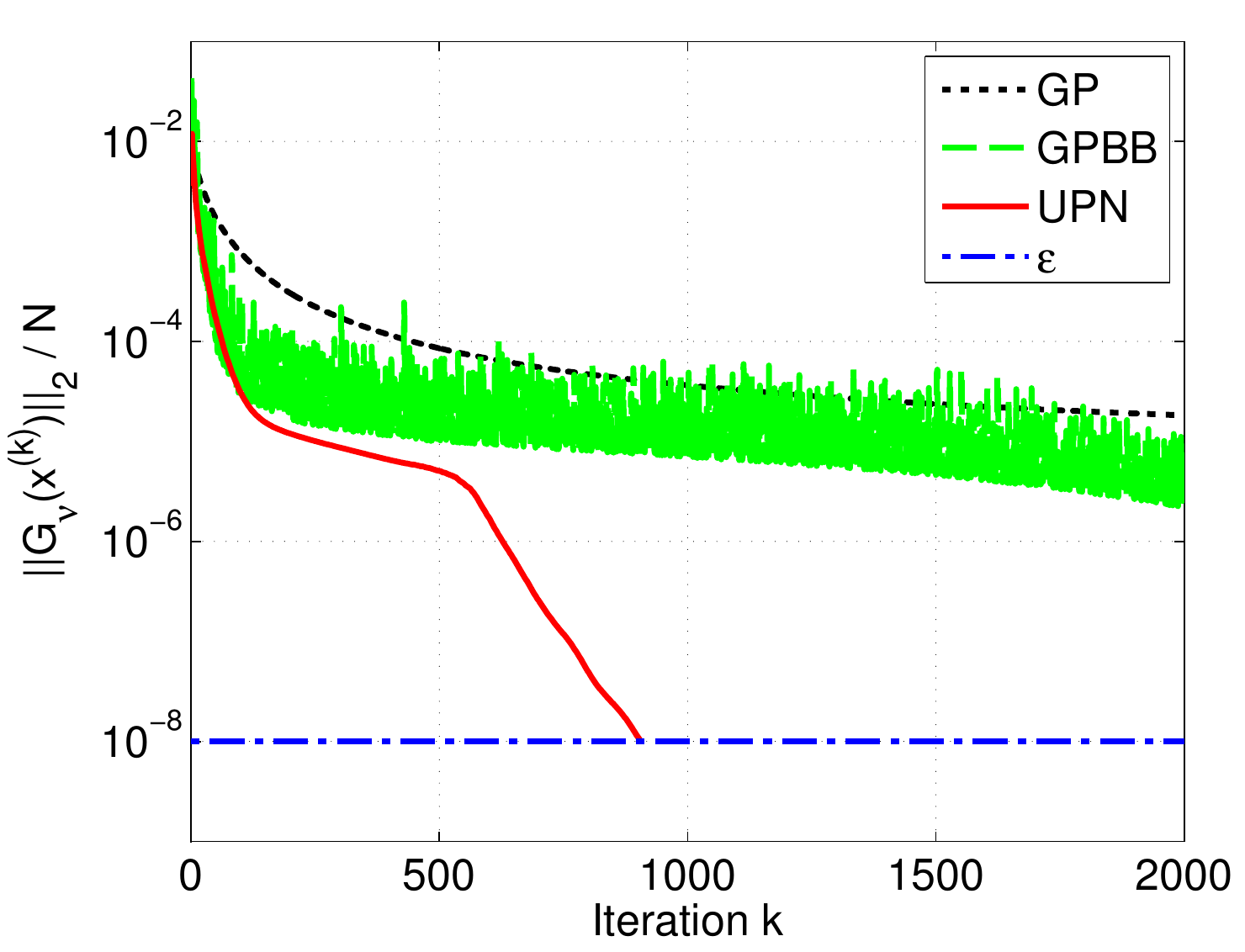}}
\caption{Convergence histories. Left: Many-view simulation. Center: Few-view simulation. Right: Gradient map norm histories for few-view simulation.}
\label{fig:convergences}
\end{figure*}

Fig. \ref{fig:reconstructions} shows the middle ($33$rd) axial voxel slice through the original 3D volume together with many-view and few-view \algname{} reconstructions using $\alpha = 0.01$. Both reconstructions reproduce the orignal features accurately, except for two small features are missing in the few-view reconstruction. 
Fig. \ref{fig:convergences} shows the convergence of the three methods in terms of objective value $\phi^{(k)}$ relative to the true minimal objective value $\phi^\star$ as function of the iterations $k$. As $\phi^\star$ is unknown, we have approximated it by computing the \algname{} solution for an $\epsilon$ two orders of magnitude lower
than the value used in the iterations.

In both cases we see that \algname{} converges to a satisfactory accuracy
within $2000$ iterations, whereas \algnameGP{} does not, and \algnameBB{} only does in the former case. In the many-view case \algname{} and  \algnameBB{} both produce a significant (and comparable) acceleration over \algnameGP{}. In the few-view case, we also observe acceleration for both, but \algname{} stands out with much faster convergence. This is in accordance with the expectation stated in Section \ref{subsec:nesterov}.

The adequacy of the stopping criterion is evaluated by a simple visual comparison of the few-view simulation gradient map norm decay (Fig. \ref{fig:convergences} right) and the objective decay (\hbox{Fig. \ref{fig:convergences} center}). Apart from the erratic decay for \algnameBB{} (which is caused by highly irregular step length selection) there is a pronounced correspondence, and we therefore consider the stopping criterion effective.
%
%

Although \algname{} was designed for strongly convex problems, we conclude that the method also works in the non-strongly convex case of having few-view data -- in fact, from the preliminary results the non-strongly convex case is where \algname{} shows its
biggest
potential by exhibiting a much faster convergence than \algnameGP{} and \algnameBB{}.

%

\section{Conclusion}
We have described the gradient-based optimization methods \algnameBB{} and \algname{} and their worst-case convergence rates. Our simulations show that both algorithms are able to significantly accelerate high-accuracy TV-based CT \hbox{image} reconstruction compared to a simple gradient method. In particular \algname{} shows much faster convergence when applied to few-view data.
The software (implementation in C with an interface to Matlab) is available from \url{http://www2.imm.dtu.dk/~pch/TVReg/}.


%
%

%

\section*{Acknowledgment}
This work is part of the project CSI: Computational \hbox{Science} in Imaging, supported by grant 274-07-0065 from the Danish Research Council for Technology and Production Sciences.



\bibliographystyle{IEEEtran.bst}
\bibliography{fully3d_references}
%
%
%

\end{document}